\documentclass[12pt]{amsart}

\numberwithin{equation}{section}

\usepackage{amssymb}
\usepackage{enumerate}

\newtheorem{thm}{Theorem}[section]
\newtheorem{lem}[thm]{Lemma}

\newtheorem{prop}[thm]{Proposition}

\def\cd{{\mathcal D}}
\def\ce{{\mathcal E}}

\def\cam{{\mathcal M}}

\def\co{{\mathcal O}}

\def\cs{{\mathcal S}}

\def\ga{{\mathfrak A}}
\def\gb{{\mathfrak B}}

\def\gd{{\mathfrak D}}

\def\gam{{\mathfrak M}}\def\gpm{{\mathfrak m}}

\def\bc{{\mathbb C}}

\def\bm{{\mathbb M}}
\def\bn{{\mathbb N}}

\def\bt{{\mathbb T}}

\def\bz{{\mathbb Z}}

\def\a{\alpha}

  \def\G{\Gamma}
\def\d{\delta}

\def\l{\lambda} \def\L{\Lambda}

\def\m{\mu}

\def\r{\rho}
\def\s{\sigma} 
\def\t{\tau}
\def\f{\varphi} \def\F{\Phi}
\def\th{\theta}  \def\Th{\Theta}
\def\om{\omega}

\newcommand{\id}{\mathop{\rm id}}

\newcommand{\tr}{\mathop{\rm Tr}}
\def\idd{{1}\!\!{\rm I}}

\newcommand{\ty}[1]{\mathop{\rm {#1}}}

\begin{document}

\title[infinite dimensional entangled Markov chains]
{infinite dimensional entangled Markov chains}
\author{Francesco Fidaleo}
\address{Francesco Fidaleo\\
Dipartimento di Matematica\\
II Universit\`{a} di Roma ``Tor Vergata''\\
Via della Ricerca Scientifica 1, 00133 Roma, Italy}
\email{{\tt fidaleo@axp.mat.uniroma2.it}}

\begin{abstract}
We continue the analysis of nontrivial examples of
quantum Markov processes. This is done by applying the construction of entangled 
Markov chains obtained from classical Markov chains with infinite 
state--space.
The formula giving the joint correlations arises from the corresponding classical formula 
by replacing the usual matrix multiplication by the Schur multiplication.
In this way, we provide nontrivial examples of entangled Markov 
chains on 
${\displaystyle\overline{\bigcup_{J\subset\bz}
\overline{\otimes}_{J}F}^{C^{*}}}$, 
$F$ being any infinite dimensional type $\ty{I}$ factor, $J$ a 
finite interval of $\bz$, and the bar the von Neumann tensor product 
between von Neumann algebras. We then have new nontrivial examples of quantum random walks
which could play a r\^ole in 
quantum information theory.

In view of applications to quantum statistical 
mechanics too,
we see that the ergodic type of an
entangled Markov chain is completely determined by the 
corresponding ergodic type of the underlying classical chain, provided 
that the latter admits an invariant probability distribution.
This result parallels the corresponding one relative to 
the finite dimensional case. 

Finally, starting from random walks on discrete ICC groups, we 
exhibit examples of quantum Markov processes based on type $\ty{II_{1}}$ 
von Neumann factors. 
\end{abstract}

\maketitle

\section{introduction}

The recent development of quantum information raised the problem 
of finding a satisfactory quantum generalization of the classical 
random walks. The relevance of this problem for quantum information 
has been emphasized in the last past years, see e.g. 
\cite{0}--\cite{2}, \cite{3}, \cite{4}, \cite{5}--\cite{9}, 
\cite{10} for different 
solutions of this problem. However, these proposals introduce some 
features which are not 
quite satisfactory from the mathematical point of view. 

First of all, these constructions are based on special models and 
a general mathematical definition of quantum random walk seems 
to be lacking. 
Second, all these constructions are based on a quantum evolution, 
unitary in some cases, irreversible in others. Now, random walks are 
particular cases of Markov processes and it is well known that, while 
in the classical case a Markov evolution uniquely determines the law 
of the corresponding stochastic process, this is in general false in 
the quantum case. 
Finally, a desirable requirement for a quantum extension of a family 
of classical processes, is that there should exist a standard procedure 
to embed the original classical family into its quantum extension.
In a satisfactory quantum generalization of the classical random walks, 
all these requirements should be precisely formulated and fulfilled.

We dealth with these problems in the previous paper 
\cite{AFi3}, by defining a nontrivial quantum lifting of classical 
Markov chains by using the Schur multiplication (known also as the 
Hadamard multiplication, see e.g. \cite{B}). In the above mentioned 
paper, we outlined this construction for 
general classical Markov chains, concentrating our attention to 
Markov chains with finite state--space.   

In the present paper we study in detail the 
construction of entangled Markov chains based on classical Markov 
chains with infinite state--space. We refer the reader to \cite{AFi3} 
and the above quoted papers for a discussion about the motivations, 
the potential applications to quantum 
information, and further details. 

We obtain nontrivial examples of random walks 
satisfying some natural requirements. Namely, they are quantum Markov 
chains such that their restrictions to at least one maximal Abelian subalgebra, 
are classical random walks. They are 
uniquely determined, up to arbitrary 
phases, by these classical restrictions. Finally, taking into account 
possible applications to information theory, they are purely generated 
(i.e. generated by isometries, see \cite{FNW1, FNW2}, and \cite{AFi3} 
for a discussion about this point).

The entangled Markov chains so constructed, generate in a natural 
way, states on the $C^{*}$--infinite tensor product   
${\displaystyle\overline{\bigcup_{J\subset\bz}
\overline{\otimes}_{J}F}^{C^{*}}}$ (denoted in the sequel
${\displaystyle\bigotimes_{\bz}F}$ by an abuse of notation), 
provided that the underlying classical 
Markov chains admit invariant distributions. Here, 
$F$ is any infinite dimensional type $\ty{I}$ factor, $J$ is a 
finite interval of $\bz$, and the bar denotes the von Neumann tensor product 
between von Neumann algebras. 

The present paper is organized as follows. Section \ref{s2} is 
devoted to prove that the Schur multiplication is well defined also 
in infinite dimensional case. Namely, the Schur multiplication
generates a (commutative) multiplication also
for infinite dimensional type 
$\ty{I}$ factors. In Section \ref{s3}, we extend to the infinite 
dimensional case, the
definition of entangled 
lifting of classical Markov chains by using the Schur 
multiplication. The ergodic properties of the entangled Markov chains 
with a stationary distribution are investigated in Section \ref{s3a}.
As for finite dimensional case, we see that the ergodic properties of the 
entangled chain are determined by those of the underlying classical 
chain. Contrary to finite dimensional case, the problem concerning the 
pureness of strongly clustering infinite dimensional entangled Markov 
chains is left open. 
Yet, some very simple examples (see Section \ref{s4}) of pure Markov chains go 
towards the conjecture that any strongly clustering entangled Markov 
chain should generate a pure state on ${\displaystyle\bigotimes_{\bz}F}$.
Section \ref{s4} contains also the description
of entangled Markov processes arising from classical random walks 
on discrete groups. Then, starting from random walks on discrete 
Infinite Conjugacy Class (ICC for short) groups, we 
provide examples of quantum random walks based on type $\ty{II_{1}}$ 
factors.

\section{the Schur multiplication in the infinite dimensional case}
\label{s2}

Let $I$ be an index set. Consider the space 
$\stackrel{\circ}{\bm}_{I}$ consisting of all the $I\times I$ matrices 
with complex entries. The subset of all bounded $I\times I$ matrices is denoted 
as $\bm_{I}$. It gives the most general type $\ty{I}$ $W^{*}$--factor, 
as $I$ varies among all cardinalities. Define the map 
$\F:\bm_{I}\mapsto\stackrel{\circ}{\bm}_{I\times I}$ as  
\begin{equation}
\label{1}
\F(A)_{(i,j)(k,l)}:=A_{ik}\d_{ij}\d_{kl}\,.
\end{equation}
Notice that the map $\F$ is identity preserving. 

Let $A,B\in\bm_{I}$. We can define the {\it Schur multiplication} as
\begin{equation*}
(A\diamond B)_{ij}:=A_{ij}B_{ij}\,.
\end{equation*}
Taking into account that
$$
(A\otimes B)_{(i,j)(k,l)}=A_{ik}B_{jl}\,,
$$
we can extend the Schur multiplication to a map 
$m:\bm_{I\times I}\mapsto \stackrel{\circ}{\bm}_{I}$ by putting
\begin{equation}
\label{13a}
m(X)_{ij}:=X_{(i,i)(j,j)}\,.
\end{equation}
\begin{lem}
\label{3}
Let $A\in\bm_{I}$. Then $\F(A)\in\bm_{I\times I}$. 
\end{lem}
\begin{proof}
We compute
\begin{align*}
|\langle\F(A)x,y\rangle|
\equiv&|\sum\overline{y_{ij}}A_{ik}\d_{ij}\d_{kl}x_{kl}|\\
=|\sum\overline{y_{ii}}A_{ik}x_{kk}|
\leq&\|A\|_{\bm_{I}}\|x\|_{\ell^{2}(I\times I)}\|y\|_{\ell^{2}(I\times I)}
\end{align*}
which leads to $\|\F(A)\|_{\bm_{I\times I}}\leq\|A\|_{\bm_{I}}$.
\end{proof}

From now on, we consider the map $\F$ as a bounded map from 
$\bm_{I}$ to $\bm_{I\times I}$. 
\begin{prop}
\label{4}
The map $\F:\bm_{I}\mapsto\bm_{I\times I}$ is an
identity preserving $*$--morphism.
\end{prop}
\begin{proof}
It is easily seen that $\F$ is $*$--preserving 
and multiplicative.
\end{proof}

Let $\r\in L^{1}(\bm_{I})_{+}$. It is immediate to show that
\begin{equation}
\label{12}
\tr\otimes\tr(\F(\r))=\tr(\r)
\end{equation}
that is the positive element $\F(\r)$ is trace--class. Furthermore, if 
$\s\in L^{1}(\bm_{I\times I})_{+}$, then 
\begin{equation}
\label{121212}
0\leq\tr(m(\s))\leq\tr\otimes\tr(\s)\,,
\end{equation}
that is the positive element $m(\s)$ is trace--class as well.

For the sake of completeness we report the following
\begin{thm}
\label{11}
Let $\F:\bm_{I}\mapsto\bm_{I\times I}$ be the linear map given in \eqref{1}.

Then $\F$ is a normal faithful identity preserving $*$--morphism 
of $\bm_{I}$ into $\bm_{I\times I}$. 

Moreover if $p\geq1$, $\F$ restricts itself to 
isometries of $L^{p}(\bm_{I})$ into $L^{p}(\bm_{I\times I})$.
\end{thm}
\begin{proof}
Let $\r\in L^{1}(\bm_{I\times I})_{+}$, and $\{X_{\a}\}\subset\bm_{I}$ 
be a net converging to $X\in\bm_{I}$
in the $*$--weak topology. We have by \eqref{121212}, 
\begin{align*}
\lim_{\a}\tr\otimes\tr(\r\F(X_{\a}))&=\lim_{\a}\tr(m(\r)X_{\a})\\
=\tr(m(\r)X)&=\tr\otimes\tr(\r\F(X))\,,
\end{align*}
where $m$ is the 
Schur multiplication. Namely, $\F$ is normal. The first part follows by Proposition 
\ref{4}, taking 
into account that a normal $*$--morphism between von Neumann factors 
is automatically faithful.

Let now $p\geq1$, and $T\in L^{p}(\bm_{I})$.
We get by the previous results,
\begin{align*}
\|\F(T)\|^{p}_{p}\equiv&\tr\otimes\tr((\F(T)^{*}\F(T))^{p/2})\\
=&\tr\otimes\tr(\F((T^{*}T)^{p/2}))\\
=&\tr((T^{*}T)^{p/2})\equiv\|T\|^{p}_{p}\,.
\end{align*}
\end{proof}

The main properties of the Schur multiplication \eqref{13a} are 
summarized in the following
\begin{thm}
The Schur multiplication \eqref{13a} defines a completely positive 
identity preserving normal map from $\bm_{I\times I}$ into $\bm_{I}$.

Moreover, $m=\big(\F\lceil_{L^{1}(\bm_{I})}\big)^{*}$ where $\F$ is the map given in 
\eqref{1}.
\end{thm}
\begin{proof}
An easy computation gives $\big(\F\lceil_{L^{1}(\bm_{I})}\big)^{*}=m$. 
Then, taking into account the 
properties of $\F$, it remains to check the normality. Let 
$\{X_{\a}\}\subset\bm_{I\times I}$ be a net converging to $X\in\bm_{I\times I}$
in the $*$--weak topology, and $\r\in L^{1}(\bm_{I})_{+}$. 
We get by \eqref{12},
\begin{align*}
\lim_{\a}\tr(\r m(X_{\a}))&=\lim_{\a}\tr\otimes\tr(\F(\r)X_{\a})\\
=\tr\otimes\tr(\F(\r)X)&=\tr(\r m(X))\,,
\end{align*}
that is $m$ is normal.
\end{proof}

Taking into account the above results, we have that
$m\lceil_{L^{1}(\bm_{I\times I})}$ is bounded and  
$\F=\big(m\lceil_{L^{1}(\bm_{I\times I})}\big)^{*}$.

From the above considerations, the Schur multiplication $m$
can be defined as a completely positive 
identity preserving normal map from $F\overline{\otimes}F$ to $F$ for arbitrary 
type I factor $F$, provided that a complete system of matrix units 
$e:=\{e_{ij}\}\subset F$ is kept fixed.

\section{infinite dimensional entangled Markov chains}
\label{s3}

Let $I$ be any index set which is kept fixed during the analysis. 
Consider a copy $\cam_{j}$ of the algebra 
$\bm_{I}\equiv\bm_{I}(\bc)$ of all bounded
$I\times I$ matrices with complex entries, together with a copy 
$\cd_{j}$ of the maximal Abelian subalgebra 
$\ell^{\infty}(I)$ of $\bm_{I}(\bc)$. For each finite subset 
$J\subset\bz$, we put  
$$
\gam_{J}:=\overline{\bigotimes}_{j\in J}\cam_{j}\,,\quad
\gd_{J}:=\overline{\bigotimes}_{j\in J}\cd_{j}\,,
$$
where $\overline{\otimes}$ is the usual von Neumann tensor product 
between von Neumann algebras. If $F\subset G$, we consider the natural 
embedding $a_F\mapsto a_F\otimes\idd_{G\backslash F}$.

The local algebra 
$$
\gam:=\overline{\big(\lim_{\stackrel{\longrightarrow}
{J\uparrow\bz}}\gam_{J}\big)}^{\,C^{*}}
$$
is the $C^{*}$--inductive limit
associated to the directed system $\{\gam_{J}\}_{J\subset\bz}$, $J$ 
finite subsets of $\bz$.
For our purpouse, we consider also the maximal Abelian subalgebra 
$$
\gd:=\overline{\big(\lim_{\stackrel{\longrightarrow}
{J\uparrow\bz}}\gd_{J}\big)}^{\,C^{*}}
$$
of $\gam$ made of the $C^{*}$--inductive limit
associated to the directed system $\{\gd_{J}\}_{J\subset\bz}$ as before.

For general $C^{*}$--algebras $\ga$, $\gb$, a
completely positive identity pre\-ser\-ving linear map 
$\ce:\ga\otimes\gb\mapsto\gb$ will be called in the sequel, a {\it transition 
expectation}.

Consider a completely positive normal map $P:\bm_{I}\mapsto\bm_{I}$. 
Such a map is said to be {\it Schur identity preserving} if 
$\ce:\bm_{I}\overline{\otimes}\bm_{I}\mapsto\bm_{I}$ given by
\begin{equation}
\label{14}
\ce=m\circ(\id\otimes P)
\end{equation}
is identity preserving 
(where ``$\circ$'' stands for composition of maps). This condition means
$$
P(\idd)_{ii}=1\,,\quad i\in I\,.
$$

Following \cite{AFi3}, Definition 1, 
any such a $P$ is said to be 
an {\it entangled Markov operator} if $P(\idd)\neq\idd$.

We 
consider also unbounded entangled Markov operators 
$P:\bm_{I}\mapsto\stackrel{\circ}{\bm}_{I}$ such that the 
associated maps 
\eqref{14} are well defined (i.e. bounded) transition expectations, see below. 

Let $P$ be an entangled Markov operator. Consider $\ce$ as above. 
Put $\ce_{A}(B):=\ce(A\otimes B)$ and 
consider a ``initial distribution'' $\r$ which is a positive normalized 
element in $L^{1}(\bm_{I})$ satisfying 
$\r=\r\circ\ce_{\idd}$.\footnote{In the infinite dimensional case, an 
invariant distribution does not always exist in general.} 
A state $\om\in\cs(\gam)$ is uniquely 
determined by all the ``finite dimensional distributions'' 
$$
\om(A_{1}\otimes\cdots\otimes A_{n}):=
\r(\ce_{A_{1}}\circ\cdots\circ\ce_{A_{n}}(\idd))\,.
$$

Such a state is a translation invariant quantum Markov chain on $\gam$
generated by the triplet $(\bm_{I},\ce,\r)$ following the terminology of 
\cite{FNW2}. 
It generalizes the construction given in \cite{AFi3} to the infinite 
dimensional case. 
For further details about the quantum Markov chains, we refer 
to \cite{A0, BrKiJoWe00, FNW1, FNW2}, 
and the references cited therein.

Now we specialize the matter to the quantum Markov chains generated in a 
canonical way by classical Markov chains with an infinite state--space, 
that is to {\it infinite dimensional 
entangled Markov chains}.  

Let $\Pi\in\stackrel{\circ}{\bm}_{I}$ be a stochastic matrix. Define, 
for $A\in\bm_{I}$, 
\begin{equation}
\label{z1}
P(A)_{ij}:=\sum_{k,l\in I}\sqrt{\Pi_{ik}\Pi_{jl}}A_{kl}\,.
\end{equation}     
By Holder inequality, we have
\begin{align*}
&\big|P(A)_{ij}\big|\equiv
\big|\sum_{k,l}\sqrt{\Pi_{ik}}A_{kl}\sqrt{\Pi_{jl}}\big|\\
&\leq\|A\|_{\bm_{I}}\big(\sum_{k}\Pi_{ik}\big)^{1/2}
\big(\sum_{k}\Pi_{jk}\big)^{1/2}\equiv\|A\|_{\bm_{I}}\,,
\end{align*}
that is $P(A)\in\stackrel{\circ}{\bm}_{I}$.\footnote{The entangled operator 
$P$ given in \eqref{z1} is not bounded 
in general. However, we define the entangled Markov chain directly in terms of the 
transition expectation given in \eqref{14}. This allows us to treat 
entangled Markov chains arising from any stochastic matrix without 
affecting our analysis, see Proposition \ref{p1}.}

We show that the transition expectation $\ce$ is purely generated, 
following  
the terminology of \cite{FNW2}.
\begin{prop}
\label{p1}
We have for the map $\ce=m\circ(\id\otimes P)$, with $P$ as 
in \eqref{z1}, 
$$
\ce(A)=V^{*}AV\,,\quad A\in\bm_{I\times I}\,,
$$
where 
the isometry $V:\ell^{2}(I)\mapsto\ell^{2}(I\times I)$ is given by
\begin{equation}
\label{rotto1}
Ve_{i}=\sum_{j\in I}\sqrt{\Pi_{ij}}e_{i}\otimes e_{j}\,,
\end{equation}
and $\{e_{i}\}_{i\in I}$ is the canonical basis of $\ell^{2}(I)$.

Hence, $\ce$ extends to a completely positive, identity preserving normal map of
$\bm_{I}\overline{\otimes}\bm_{I}\cong \bm_{I\times I}$ into $\bm_{I}$.
\end{prop}
\begin{proof}
It is immediate to show that Formula \eqref{rotto1} defines a bounded operator.
The proof follows as 
$$
V^{*}e_{i}\otimes e_{j}=\sqrt{\Pi_{ij}}e_{i}\,.
$$
\end{proof}

Let now $\pi:=\{\pi_{j}\}_{j\in I}$ be an invariant measure for $\Pi$. Define
the matrix $Q(\pi)\in\stackrel{\circ}{\bm}_{I}$ given by
\begin{equation}
\label{z2}
Q(\pi)_{ij}:=\sum_{k\in I}\pi_{k}\sqrt{\Pi_{ki}\Pi_{kj}}\,.
\end{equation}     
By Holder inequality, we get
\begin{align*}
&Q(\pi)_{ij}\equiv\sum_{k}\sqrt{\pi_{k}\Pi_{ki}}\sqrt{\pi_{k}\Pi_{kj}}\\
\leq&\big(\sum_{k}\pi_{k}\Pi_{ki}\big)^{1/2}
\big(\sum_{k}\pi_{k}\Pi_{kj}\big)^{1/2}\equiv\sqrt{\pi_{i}\pi_{j}}\,,
\end{align*}
that is, $Q(\pi)$ is well defined. Moreover, it is 
positive by construction. Furthermore, $Q(\pi)$ defines a 
bounded positive form on $\bm_{I}$ if and only if $\pi$ defines a 
positive form on $\ell^{\infty}(I)$, with 
$\|Q(\pi)\|_{L^{1}(\bm_{I})}=\|\pi\|_{\ell^{1}(I)}$.
\begin{prop}
\label{inv}
$Q$ given in \eqref{z2} maps the set of 
invariant measures 
for $\Pi$ into the set of normal semifinite weights on $\bm_{I}$ invariant for 
$\ce_{\idd}\equiv \idd\diamond P(\,\cdot\,)$, 
with $P$ given in \eqref{z1}. 

$Q$ restricts itself to 
a one--to--one correspondence between the invariant probability 
distributions 
for $\Pi$ and the normal states invariant for 
$\ce_{\idd}$.
\end{prop}
\begin{proof} 
Taking into account the proof of Proposition 2.4 of \cite{AFi3}, it 
is enough to prove semifiniteness and normality. Let $\pi$ be as 
above. Let $\gpm_{Q(\pi)}$ be the definition--domain of the weight 
$Q(\pi)$ (\cite{T2}, Definition VII.1.3). Then  
$$
\gpm_{0}:=\bigcup\big\{\bm_{F}\,\big|\,F\,\text{finite subset of}\, 
I\big\}\subset\gpm_{Q(\pi)}
$$ 
and 
$\gpm_{0}''=\bm_{I}$, that is $Q(\pi)$ is semifinite.

We write, with an abuse of notation,
$Q(\pi)(A)=\pi(\ce(\idd\otimes A))$ as $\ce(\idd\otimes A)\in\ell^{\infty}(I)$. 
The normality follows by Fatou Lemma (i.e. $\pi$ is normal on 
$\ell^{\infty}(I)$), taking into account that $\ce$ is a normal map.
\end{proof}

Let $(\Pi,\pi)$ consist of a stochastic matrix as above, and an invariant 
probability measure for it, respectively. Put $\bm:=\bm_{I}$.
The {\it entangled Markov chain} associated to
$(\Pi,\pi)$ is the translation invariant locally normal state on $\gam$ 
generated by the triplet $(\bm,\ce,Q(\pi))$. It is immediate to show 
that the quantum chain reduces itself to the classical one, when 
restricted to $\gd$.

An invariant distribution (i.e. an invariant 
probability measure) does not always exist for a given infinite 
stochastic matrix, see e.g. Theorem 6.2.1 of \cite{Ci}. However, by 
Proposition \ref{inv}, one can 
define on $\gam$ a translation invariant locally normal weight 
starting from the triplet $(\bm,\ce,Q(\pi))$, $\pi$ being an 
invariant measure for $\Pi$ if the last exists.\footnote{An invariant measure
for the stochastic matrix $\Pi$ always exists if the latter contains 
recurrent indices (``states'' in the terminology of classical Markov 
chains), see e.g. \cite{Ci}, Theorem 6.2.25 and Theorem 5.3.14.}

Notice that, by Proposition \ref{p1}, an entangled Markov chain is always 
generated by an isometry (purely generated in the terminology of 
\cite{FNW2}). 

As it was shown in \cite{AFi3}, the entangled Markov chains are 
defined up to arbitrary phases. Namely, let 
$\stackrel{\circ}{\bm}_I\!\!(\bt)$ be the set consisting of all the 
$I\times I$ matrices with entries in  
the unit circle $\bt$. If $\chi\in\bm_I\!(\bt)$, 
then the map   
$P_{\chi}$ defined as
\begin{equation}
\label{circle}
P_{\chi}(A)_{ij}:=\sum_{k,l\in I}\overline{\chi_{ik}}\chi_{jl}
\sqrt{\Pi_{ik}\Pi_{jl}}a_{kl}
\end{equation}
gives rise to an entangled Markov operator as well. The corresponding 
quantum measures $Q_{\chi}(\pi)$ associated to the entangled Markov 
operator in \eqref{circle} are given by
$$
Q_{\chi}(\pi)_{ij}=\sum_{k\in I}\pi_{k}
\chi_{ki}\overline{\chi_{kj}}\sqrt{\Pi_{ki}\Pi_{kj}}\,,
$$
provided that the classical chain $\Pi$ admits the invariant measure $\pi$.

\section{ergodic properties}
\label{s3a}

The present section is devoted to the investigation of 
the ergodic properties of infinite dimensional entangled 
Markov chains. The following analysis parallels the corresponding one 
of Section 3 of 
\cite{AFi3} relative to the finite dimensional case.

Taking into account Formula (3.15) of \cite{Ci}, we write
\begin{equation}
\label{111}
\Pi=p_0\Pi p_0+\sum_{\th\in\Th}(p_0\Pi p_\th+p_\th\Pi 
p_\th)+\sum_{\l\in\L}(p_0\Pi p_\l+p_\l\Pi p_\l)\,.
\end {equation} 

Here, $\Th$, $\L$ label the recurrent--null, and the 
recurrent--positive classes of $\Pi$, respectively, and $p_{0}$ is 
the selfadjoint projection associated to the indices relative to the 
transient states of the classical chain under consideration.  
Furthermore, for each $\Pi_{\l}:=p_\l\Pi p_\l$,
\begin{equation*}
\Pi_{\l}=\sum_{j_{\l}=1}^{m_{\l}}
p_{\l,j_{\l}}\Pi_{\l}p_{\l,j_{\l}+1}
=:\sum_{j_{\l}=1}^{m_{\l}}
p_{\l,j_{\l}}\Pi_{j_{\l},j_{\l}+1}p_{\l,j_{\l}+1}
\end{equation*}     
with the convention that $m_{\l}+1=1$. Here, $m_{\l}$ is the (finite) period 
corresponding to indices (``states'' in the terminology of classical 
Markov chains) of the ergodic class $\l$, with the convention that 
$m_{\l}=1$ in the aperiodic case. Notice that the indices (i.e. 
``states'')  
corresponding to an ergodic class are at most denumerable.
We suppose that $\L$ is nonvoid.\footnote{This is always the case in 
the finite dimensional case, where the class of recurrent--null 
states is void by finiteness, see \cite{S} for further details.} In 
this case, there exist stationary distributions for $\Pi$, and by 
Proposition \ref{inv}, stationary distributions for $\ce$ given by \eqref{14}. 
Let $\pi$ be 
any such a distribution. It has the form 
\begin{equation*}
x=\sum_{\l\in\L}\a_{\l}x_{\l}\,,\quad 
\sum_{\l\in\L}\a_{\l}=1\,,\,\, \a_{\l}\geq0\,, 
\,\,\l\in\L\,,
\end{equation*}
$x_{\l}$ being the unique stationary distribution for the stochastic 
matrix $\Pi_{\l}$.
Let $\L_{\pi}\subset\L$ be the set of the ergodic classes $\l$ of $\Pi$ 
such that
$\a_{\l}>0$. Of 
course, the cardinality of $\L_{\pi}$ is at most denumerable. If
\begin{equation}
\label{1113}
p:=\sum_{\l\in\L}p_{\l}
\end{equation}
is the support in $\ell^{\infty}(I)$ of $\pi$ considered in a natural 
way as an element of $\bm$,
define $\tilde\ce:\bm\otimes\bm_{p}\mapsto\bm_{p}$
given by
\begin{equation}
\label{113}
\tilde\ce:=p\ce\lceil_{\bm_{p}}(\,\cdot\,)p\,.
\end{equation}     

Let $\eta$ be the action of the completely reducible part $p\Pi p$ on 
the set of projections 
$\big\{\{p_{\l,j_{1}},
\dots,p_{\l,j_{m_{\l}}}\}\big\}_{\l\in\L}$, $p$ being given in 
\eqref{1113}.
Such an action leaves each ergodic component 
$\{p_{\l,j_{\l}}\}_{\l=1}^{m_{\l}}$ globally invariant, acting 
cyclically on it. Choose any projection, say 
$\bar p_{\l}:=p_{\l,\bar j}$, in each ergodic class $\l\in\L_{\pi}$,
where $\L_{\pi}$ labels the ergodic components present in the stationary 
distribution $\pi$ as described above. Define for $\{A_{1},\dots,A_{n}\}\subset\bm$,
\begin{align*}
&\f_{\l}(A_{1}\otimes\cdots\otimes A_{n}):=\pi(\bar p_{\l})^{-1}
Q(\pi)(\ce_{A_{1}}\circ\cdots\circ\ce_{A_{n}}(\eta^{-n}\bar p_{\l}))\,,\\
&\om_{\l}(A_{1}\otimes\cdots\otimes A_{n}):=\pi(p_{\l})^{-1}
Q(\pi)(\ce_{A_{1}}\circ\cdots\circ\ce_{A_{n}}(p_{\l}))\,.
\end{align*}

Let $\om$ be the entangled Markov chain on $\gam$ associated to the 
triplet $(\bm,\ce,Q(\pi))$. 
It is straightforward to verify that
$$
\om_{\l}=\frac{1}{m_{\l}}\sum_{k=1}^{m_{\l}}
\f_{\l}\circ\t^{k}\,,
$$
\begin{equation}
\label{cczz}
\om=\sum_{\l\in\L_{\pi}}\frac{\pi(p_{\l})}{m_{\l}}\sum_{k=1}^{m_{\l}}
\f_{\l}\circ\t^{k}\,,
\end{equation} 
where $\t$ is the one--step shift on the chain.

The states $\om_{\l}$, $\f_{\l}$ describe the decomposition of $\om$ into 
ergodic and completely ergodic components, 
respectively.\footnote{The state $\f_{\l}$ is only 
$m_{\l}$--step translation invariant, $m_{\l}$ being the period of 
the component $\l$, and keeps track of the localization (modulo a 
period), see \cite{AFi3}, Section 5 for the precise way to define 
$\f_{\l}$ on $\gam$.} 
\begin{thm}
\label{main}
Let $(\Pi,\pi)$ consist of a stochastic matrix 
and a stationary distribution for it, with $\L_{\pi}\neq\emptyset$. 
Consider the entangled Markov chain 
$\om$ on $\gam$ generated by the triplet $(\bm,\ce,Q(\pi))$. 
The following assertions hold true.
\begin{itemize}
\item[(i)] The state $\om$ is ergodic w.r.t the spatial translations 
if and only if the set $\L_{\pi}$ is a singleton. 
\item[(ii)] The state $\om$ is strongly clustering w.r.t the spatial translations
if and only if the set $\L_{\pi}$ is a singleton, and in addition, 
the corresponding block in 
the decomposition \eqref{111} of $\Pi$ is aperiodic. 
\end{itemize}
\end{thm}
\begin{proof}
It is immediate to verify that
$\om$ is given by \eqref{cczz}. Furthermore, the states appearing in the 
r.h.s. of \eqref{cczz} give rise to different states when restricted to the 
Abelian subalgebra $\gd$
of $\gam$. So, they are mutually different. It is then enough to show 
that the $\om_{\l}$ are ergodic w.r.t. the one--step shift, and
the $\f_{\l}$ are strongly clustering w.r.t. the $m_{\l}$--step shift, 
respectively. 

Let $A=A_{1}\otimes\cdots\otimes A_{r}$, $B=B_{1}\otimes\cdots\otimes B_{s}$, 
we compute by applying Lemma 2.1 of \cite{AFi3},
\begin{align*}
&\frac{1}{n}\sum_{k=1}^{n}\om_{\l}(A\t^{k}(B))
=\pi(p_{\l})^{-1}\\
\times&Q(\pi)\bigg(\ce_{A_{1}}\circ\cdots\circ\ce_{A_{r}}
\circ\bigg(\frac{1}{n}\sum_{k=0}^{n-1}\ce_{\idd}^{k}\bigg)
\bigg(\ce_{\idd}\circ\tilde\ce_{B_{1}}\circ\cdots\circ\tilde\ce_{B_{s}}
(p_{\l})\bigg)\bigg)\,,
\end{align*}
where $\tilde\ce$ is given by \eqref{113}. Define, componentwise, the element 
$v\in\ell^{1}(I)$ as 
$$
v_{i}:=Q(\pi)\big(\ce_{A_{1}}\circ\cdots\circ\ce_{A_{r}}(e_{ii})\big)\,.
$$
Define the element $D\in\ell^{\infty}(I)_{p_{\l}}\subset\bm_{I}$ as
$$
D:=\ce_{\idd}\circ\tilde\ce_{B_{1}}\circ\cdots\circ\tilde\ce_{B_{s}}(p_{\l})\,.
$$
We get by Dominated 
Convergence Theorem,  
\begin{align*}
&\lim_{n}Q(\pi)\bigg(\ce_{A_{1}}\circ\cdots\circ\ce_{A_{r}}
\circ\bigg(\frac{1}{n}\sum_{k=0}^{n-1}\ce_{\idd}^{k}(D)\bigg)\bigg)\\
=&\lim_{n}Q(\pi)\bigg(\ce_{A_{1}}\circ\cdots\circ\ce_{A_{r}}
\circ\bigg(\frac{1}{n}\sum_{k=0}^{n-1}\Pi_{\l}^{k}D\bigg)\bigg)\\
=&\lim_{n}\bigg\langle v,\bigg(\frac{1}{n}
\sum_{k=0}^{n-1}\Pi_{\l}^{k}\bigg)D\bigg\rangle\\
=&\pi(p_{\l})^{-1}Q(\pi)(\ce_{B_{1}}\circ\cdots\circ\ce_{B_{s}}
(p_{\l}))
Q(\pi)(\ce_{A_{1}}\circ\cdots\circ\ce_{A_{r}}(p_{\l}))\,.
\end{align*}
Here, $\langle\,\cdot\,,\,\cdot\,\rangle$ is the natural pairing 
between $\ell^{1}$ and $\ell^{\infty}$, and 
the last equality follows by Theorem 6.1 of \cite{KSK}.

Collecting together, we have
$$
\lim_{n}\frac{1}{n}\sum_{k=1}^{n}\om_{\l}(A\t^{k}(B))=\om_{\l}(A)\om_{\l}(B)\,,
$$
that is, $\om_{\l}$ is ergodic. 
The mixing property w.r.t. the $m_{\l}$--step shift for the $\f_{\l}$ 
is proven in the same way, taking into account Theorem 6.38 of \cite{KSK}.

Now, the quantum chain is ergodic iff $\L_{\pi}$ is a singleton (i.e. 
$\L_{\pi}=\{\l_{0}\}$),  which 
corresponds to the case when there is only one summand in the first 
sum in \eqref{cczz}.
The quantum chain is strongly clustering when, in addition, also the second sum in
\eqref{cczz} consists of one element, that is when, for the period, $m_{\l_{0}}=1$.  
\end{proof}

\section{some applications}
\label{s4}

We are going to consider some interesting examples of 
infinite dimensional entangled Markov chains.

\medskip

We begin with the classical chain $(\Pi,\pi)$, $\pi$ being an 
invariant distribution for the stochastic matrix $\Pi$.
Consider the collections 
$\{D^{(k)}\}_{k\in\bn\backslash\{0\}}$ of 
the density matrices relative to the local algebras
$$
\gam_{[1,k]}:=\overline{\bigotimes}_{1\leq j\leq k}\cam_{j}
$$
and their translates, 
arising from the entangled Markov chain $(\bm,\ce,Q(\pi))$ associated 
to $(\Pi,\pi)$.
Here, $\ce$ is given in \eqref{14}, with the 
entangled operator $P$ given in \eqref{z1}. 
 
\begin{prop}
\label{dms}
We have for the collections 
$\{D^{(k)}\}_{k\in\bn\backslash\{0\}}$ of 
the density matrices,
\begin{equation}
\label{elm}
D^{(k)}_{(i_{1},\dots,i_{k})(j_{1},\dots,j_{k})}=
\tr{}_{\bm}\bigg(\ce_{D_{\pi}}\circ\ce_{e_{i_{1}j_{1}}}
\circ\cdots\circ\ce_{e_{i_{k}j_{k}}}(\idd)\bigg)\,,
\end{equation}
where $D_{\pi}:=\sum_{k\in I}\pi_{k}e_{kk}$ is the diagonal embedding 
of $\pi$ in $L^{1}(\bm)$, and the $e_{ij}$ are the canonical matrix units 
of $\bm$.
\end{prop}
\begin{proof}
A simple computation.
\end{proof}

We specialize the situation to the entangled Markov chain generated 
by the the aperiodic stochastic projection $Q$ with
matrix elements $q_{ij}=\pi_{j}>0$, $i\in\bn$ and
$\pi\in\ell^{1}$. We obtain for the sequences of density matrices given in 
\eqref{elm},
$$
D^{(k)}=\underbrace{D^{(1)}\otimes\cdots\otimes D^{(1)}}_{k\text{--times}}\,,
$$
with $D^{(1)}_{ij}=\sqrt{\pi_{i}\pi_{j}}$. The state $\om$ 
on $\gam$ associated to these examples of entangled Markov chains, is then an infinite 
product vector state based on the vector 
$\{\sqrt{\pi_{i}}\}_{i\in\bn}\in\ell^{2}$ (\cite{vN}). 
Namely, $\om$ is a pure state, see \cite{Sa}, Proposition 4.4.4.\footnote{This non generic 
situation exhibits a low degree of entanglement. These examples are 
then not 
suitable for possible applications to quantum information theory.}

This very simple example, together with the result in Theorem 3.4 of 
\cite{AFi3} relative to the pureness of finite dimensional entangled 
Markov chains, allows us to conjecture that strongly 
clustering infinite dimensional entangled Markov chains generate pure 
states on $\gam$. This might be proved by studying the algebraic properties of 
the sequence of the ranges of the $D^{(k)}$ (see \cite{FNW1}), 
taking into account 
\eqref{elm}. For this aim, we report another useful formula for the density 
matrices $D^{(k)}$. 
Define the matrix $\G\in\stackrel{\circ}{\bm}_{I\times I}$ as
$$
\G_{(i,j)(k,l)}:=\sqrt{\Pi_{ij}\Pi_{kl}}\,.
$$
We get
$$
D^{(k)}=\underbrace{m\otimes\cdots\otimes m}_{k\text{--times}}
\big(Q(\pi)\otimes\underbrace{\G\otimes\cdots\otimes 
\G}_{(k-1)\text{--times}}\otimes P(\idd)\big)
$$
where $m$ is the Schur multiplication, $P$ is the entangled Markov 
operator associated to the stochastic matrix $\Pi$ via \eqref{z1}, and $Q(\pi)$
is the trace--class matrix associated to the stationary distribution 
$\pi$ via \eqref{z2}.

\medskip

We pass to the entangled Markov processes based on random walks on 
discrete groups. We note that most of such Markov chains  
generate merely locally normal weights on $\gam$. For some basic facts 
about random walks on groups, see e.g. \cite{R}.

Let $G$ be a discrete group and $\m$ a probability measure on it 
which is kept fixed during the analysis. The {\it right} and 
{\it left} random walks on $G$ are given by the doubly stochastic
transition matrices  
$\Pi^{r}$, $\Pi^{l}$ respectively, with
$$
\Pi^{r}_{gh}:=\m(g^{-1}h)\,,\quad \Pi^{l}_{gh}:=\m(gh^{-1})\,.
$$

Let $P^{r}$, $P^{l}$ be the corresponding entangled Markov operators
obtained by \eqref{z1}.\footnote{As $\Pi^{r}$, $\Pi^{l}$ are doubly 
stochastic, the associated entangled operators $P^{r}$, $P^{l}$ are 
bounded.} An easy computation gives rise for $g\in G$,
\begin{equation}
\label{tras}
\text{ad}(\l(g))\circ P^{r}=P^{r}\circ\text{ad}(\l(g))\,,\, 
\text{ad}(\r(g))\circ P^{l}=P^{l}\circ\text{ad}(\r(g))\,,
\end{equation}
where $\l$, $\r$ are the left and right translations on $G$. Here, 
\eqref{tras} follows by the corresponding equivariance properties of 
$\Pi^{r}$, $\Pi^{l}$. Denote by $R(G)$, $L(G)$ the von Neumann 
algebras generated by the right and left translations on $G$.
\begin{prop}
We have for the entangled Markov operators $P^{r}$, $P^{l}$, and 
for the transition expectations $\ce^{r}$, $\ce^{l}$ constructed as 
in \eqref{14},
\begin{align*}
&P^{r}(R(G))\subset R(G)\,,\quad P^{l}(L(G))\subset L(G)\,,\\  
&\ce^{r}(R(G)\overline{\otimes}R(G))\subset R(G)\,,\quad 
\ce^{l}(L(G)\overline{\otimes}L(G))\subset L(G)\,.
\end{align*}
\end{prop}
\begin{proof}
Taking into account that $\r(g)_{xy}=\d_{x,yg^{-1}}$, 
$\l(g)_{xy}=\d_{x,gy}$, we easily obtain
$$
P^{r}(R(G))\subset L(G)'\equiv R(G)\,,
\quad P^{l}(L(G))\subset R(G)'\equiv L(G)\,.
$$
The proof follows as, for the Schur multiplication,
$$
m(R(G)\overline{\otimes}R(G))\subset R(G)\,,\quad 
m(L(G)\overline{\otimes}L(G))\subset L(G)\,.
$$
\end{proof}

We end by noticing that, starting with a discrete ICC group, one 
can construct entangled Markov processes on 
${\displaystyle\overline{\bigcup_{J\subset\bz}
\overline{\otimes}_{J}R}^{C^{*}}}$ ($J$ runs on all finite subsets of $\bz$), 
$R$ being the algebra $L(G)$ 
generated by left translations (equally well the algebra $R(G)$ 
generated by right translations). Then, we provide entangled Markov 
processes based on type $\ty{II_{1}}$ factors. 

We hope to return elsewhere on all the questions left open in this 
section. 

\section*{acknowledgements}

The author is grateful to L. Accardi for suggesting the problem and
for fruitful discussions.

\end{document}